\documentclass[11]{siamltex}
\usepackage{amsmath,amssymb,amsfonts,graphicx,fancybox,shadow,color}

\newtheorem{alg}{Algorithm}

\usepackage{graphicx}

\def\be#1\ee{\begin{equation}#1\end{equation}}

\begin{document}
\bibliographystyle{plain}

\pagestyle{myheadings}

\title{ROM inversion of monostatic data lifted to full MIMO. }

\author{
 V. Druskin\footnotemark[1],  S. Moskow\footnotemark[2] and M. Zaslavsky\footnotemark[3]}

\renewcommand{\thefootnote}{\fnsymbol{footnote}}

\footnotetext[1]{Worcester Polytechnic Institute, Department of Mathematical Sciences,
Stratton Hall,
100 Institute Road, Worcester MA, 01609 (vdruskin1@gmail.com)}
\footnotetext[2]{Department of Mathematics, Drexel University, Korman Center, 3141 Chestnut Street, Philadelphia, PA 19104
(moskow@math.drexel.edu)}
\footnotetext[3]{Southern Methodist University, Department of Mathematics, Clements Hall, 3100 Dyer st,
Dallas, TX 75205 (mzaslavskiy@smu.edu)}
\maketitle

\begin{abstract} The Lippmann--Schwinger--Lanczos (LSL) algorithm has recently been shown to provide an efficient tool for imaging and direct inversion of synthetic aperture radar data in multi-scattering environments \cite{DrMoZa3}, where the data set is limited to the monostatic, a.k.a. single input/single output (SISO) measurements. The approach is based on constructing data-driven estimates of internal fields via a reduced-order model (ROM) framework and then plugging them into the Lippmann-Schwinger integral equation. However, the approximations of the internal solutions may have more error due to missing the off diagonal elements of the multiple input/multiple output (MIMO) matrix valued transfer function. This, in turn, may result in multiple echoes in the image. Here we present a ROM-based data completion algorithm to mitigate this problem. First, we apply the LSL algorithm to the SISO data as in \cite{DrMoZa3} to obtain approximate reconstructions as well as the estimate of internal field. Next, we use these estimates to calculate a forward Lippmann-Schwinger integral to populate the missing off-diagonal data (the lifting step). Finally, to update the reconstructions, we solve the Lippmann-Schwinger equation using the original SISO data, where the internal fields are constructed from the lifted MIMO data. The steps of obtaining the approximate reconstructions and internal fields and populating the missing MIMO data entries can be repeated for complex models to improve the images even further. Efficiency of the proposed approach is demonstrated on 2D and 2.5D  numerical examples, where we see reconstructions are improved substantially. 

 \end{abstract}

\section{Introduction}
Although the Born approximation is a method of choice for synthetic aperture radar (SAR) due to its speed and versatility,  nonlinear scattering is known to produce significant artifacts.
Approaches to deal with multiple scattering include optimization \cite{9357476},  improved linearizations using some knowledge of an inhomogeneous background \cite{gilman2015mathematical,WaAnLi}, linear sampling methods \cite{CaCoMo, Burfeindt2}, statistical and machine learning methods \cite{lagergren2021deep}, iterative time reversal and matching filtering \cite{10352835}, and reweighting \cite{cheney2021procedure}.
Also in seismic imaging are  full wave inversion \cite{Virieux2016AnIT}, boundary control \cite{Kepley2016GeneratingVI}, and inverse series methods \cite{Malcolm2007IdentificationOI}.
Many of the above approaches require prior knowledge, the solution of full scale forward problems or large training sets.
%

Reduced order model (ROM) approaches have been shown to provide accurate images for several inverse impedance, scattering and diffusion problems \cite{borcea2011resistor,borcea2014model,druskin2016direct,druskin2018nonlinear,borcea2017untangling,borcea2019robust,BoDrMaMoZa,
borcea2020reduced, DrMoZa,  Borcea2021ReducedOM, DrMoZa2, BoGaMaZi}.  In the recent work \cite{DrMoZa3}, we applied  the Lippmann-Schwinger-Lanczos algorithm  \cite{DrMoZa}  to models of synthetic aperture radar.  
Most of the ROM-based imaging approaches strongly rely on having a complete square MIMO dataset accessible,  as seen in the schematic \eqref{fullMIMO}
	\begin{equation}\label{fullMIMO}\begin{pmatrix}.& . & . & . & . \\
		\  & . & . & . & . \\
		\  &  \  & . & . & . \\
		\  &  \  &  \  & . & . \\
		\  &  \  &  \  &  \  & . \\ 
		\end{pmatrix}  \end{equation}
where the column and row numbers correspond respectively to the indices of the receivers (outputs) and transmitters (inputs). The other half of the transfer matrix extends by reciprocity. In SAR, as well as in conventional (non-array) medical ultrasound, profiling is done by mono-static transmitter/receivers, and images are constructed by combining multiple single-input/single output (SISO) data. The data in this case is visualized in \eqref{SAR},
	\begin{equation}\label{SAR}\begin{pmatrix}.&   &   &   &   \\
		\  & . &   &   &   \\
		\  &  \  & . &   &  \\
		\  &  \  &  \  & . &   \\
		\  &  \  &  \  &  \  & . \\ 
		\end{pmatrix} \end{equation}
where we have only the diagonal. The nature of the Lippmann-Schwinger-Lanczos (LSL) algorithm allowed us to extend the ROM approach naturally to this case in \cite{DrMoZa3}.  LSL for monostatic data \cite{DrMoZa3} contains two major steps: \begin{itemize} \item{Step (i)}  Use the individual SISO ROMs to generate internal fields \item{Step (ii)} Use these internal fields to linearize the Lippmann-Schwinger equation and solve for the profile.  \end{itemize} Using this method we were able to substantially reduce the echos compared to the Born approximation, with minimal increase to the computational cost. However, due to the lack of the off diagonal data, we did not have the full strength of the MIMO ROM,  and so the data generated internal fields were less accurate. In this work we propose to use the SISO ROMs to lift the data to improve the data generated internal fields. We start by first implementing LSL for monostatic data. Once we have approximate internal fields and profile, we do a {\it forward} application of the Lippmann-Schwinger integral to approximate the off-diagonals. That is, we complete the data so that it is of the form (\ref{SARlifted})
\begin{equation}\label{SARlifted}\begin{pmatrix}.& \circ  &  \circ & \circ   & \circ  \\
		\circ  & . &  \circ & \circ  & \circ  \\
		\circ &  \circ & . & \circ  & \circ \\
		\circ  &  \circ  &  \circ  & . &  \circ \\
		\circ  &  \circ  &  \circ  &  \circ  & . \\ 
		\end{pmatrix} \end{equation}
where the symbol $\circ$ represents approximate, or lifted, data. With a full (yet approximate) transfer matrix in hand, we can use a MIMO ROM to construct improved internal fields, for a new Step (i).  These improved fields are then used in the Lippmann Schwinger equation in Step (ii). Evidently, at Step (ii), we use only the true diagonal data, since the lifted off diagonals will not provide any new information. We note that a related procedure has been recently developed for a so-called monostatic-to-multistatic transform with application to linear sampling-based imaging of SAR data \cite{Burfeindt2}.

We note that the lifted data matrix (\ref{SARlifted}) does not have the correct structure apriori,  that is to say, the resulting mass matrix may not be symmetric positive definite. Regularization will therefore be necessary at this step separately, in addition to the usual regularization required for the inversion of the Lippmann-Schwinger equation. 
We also emphasize that lifting doesn't provide additional data for the inverse problem; it is used solely to improve the ROM's performance in generating the internal fields. 

{We now make a remark about ROM based inversion methods. In general, a complete ROM is a reduced version of the full scale PDE in the form of a Galerkin system, containing both a mass matrix $M$ (the Grammian used here) and a stiffness matrix. Historically, the complete reduced order model was needed to do inversion (see for example \cite{druskin2016direct}). It was later found that in the time domain the stiffness matrix was not needed for the LSL approach, hence in this work we discuss only $M$, and refer to this as the ROM.}

This paper is organized as follows. In Section 2 we state precisely the model and inverse problem and give an overview of the entire process. Sections  3,4, and 5 contain detailed descriptions of the steps, which are monostatic LSL, data completion, and MIMO generation of internal fields, respectively. In Section 6 we discuss iteration, and in Section 7 we present numerical experiments demonstrating the improved performance. We present a brief conclusion in Section 8. As a notational note, throughout the manuscript,  we use boldface to denote data generated approximations while the regular script denotes actual values.

\section{Statement of the inverse problem and algorithm overview}
We assume the following wave model problem for a domain $\Omega\subset\mathbb{R}^n$
\be\label{waveeq}
u_{tt}  + Au=   0 \ \mbox{in}  \ \Omega\times [ 0,\infty)
\ee
with initial conditions
\begin{eqnarray} \label{waveeqic}
u ( t=0) &=& g_i \ \mbox{in}\  \Omega  \\
u_t ( t=0) &=&  0\  \mbox{in}\  \Omega \end{eqnarray}
where   $\{ g_i\} $ for $i=1,\ldots K$ are Dirichlet initial conditions representing localized sources near the boundary,\footnote{The formulation for the homogeneous wave equation \eqref{waveeq} with an inhomogeneous initial condition \eqref{waveeqic} can be equivalently derived from the more conventional  radar formulation of an inhomogeneous wave equation with homogeneous initial conditions.   We refer to \cite{druskin2016direct} for details. }    and we assume homogeneous Neumann boundary conditions on the spatial boundary $\partial \Omega$. {The operator $A$ is taken to be the negative Laplacian  plus a scattering term: \begin{equation}\label{Aeq}  A = -\Delta + q \end{equation}
where $q(x)\ge 0$ is our unknown potential, not necessarily small, but with compact support.  The initial data is assumed to be localized at source $i$, and the corresponding exact forward solutions to (\ref{waveeq}-\ref{waveeqic}) are
\begin{equation}\label{exactsolution}  u^{(i)}(x,t) = \cos{(\sqrt{A}t)} g_i(x), \end{equation}
where the square root and cosine are defined via the spectral theorem.} These solutions are assumed to be unknown, except near the receiver.  The data here is monostatic, that is, we measure the data from source $i$ only back at the receiver collocated with the source. The data is measured at the $2n-1$ evenly spaced time steps $t= k\tau$ for $k=0,\dots, 2n-2$ and is modeled by 
\begin{eqnarray}\label{datadef} F^{ii}(k\tau ) &= &\int_\Omega g_i(x) u^{(i)}(x,k\tau) dx\nonumber \\ &=&   \int_\Omega g_i(x) \cos{(\sqrt{A}k\tau)} g_i(x) dx .\end{eqnarray}
The inverse problem is as follows:  Given $$ \{ F^{ii}(k \tau) \} \  \mbox{ for } \  k=0,\dots, 2n-2; \ \  i=1,\ldots K,$$
reconstruct $q$. 
We will also need to assume we have background fields $\{ u^{(i),0} \}$, solutions to (\ref{waveeq}-\ref{waveeqic}) corresponding to $q=0$. Denote the full $K\times K$ background transfer function by $F_0(t)$ and background solution anti-derivatives by  $$w^{(i),0}(x,t)=\int^t_0{u^{(i),0}(x,s)ds}.$$  Note that $\{ w^{(i),0} \}$ satisfy (\ref{waveeq}) for $q=0$ with homogeneous Neumann boundary conditions and Neumann initial conditions
\begin{eqnarray} \label{waveeqic2}
w^{(i),0} ( t=0) &=& 0\ \mbox{in}\  \Omega  \\
w^{(i),0}_t ( t=0) &=&  g_i\  \mbox{in}\  \Omega \end{eqnarray}

Throughout the process we will use the following form of the Lippmann-Schwinger equation 
\begin{equation} \label{eq:LSTD} F^{ij}_{0}(t) - F^{ij}(t) = \int_0^{t}\int_\Omega w^{(j),0}(x,t-\tau) u^{(i)}(x,\tau) q(x) dxd\tau. \end{equation}
We remind the reader that the inverse Born method is to replace the unknown field $u^{(i)}$ in (\ref{eq:LSTD}) with the corresponding background solution $u^{(i),0}$, and solve for $q$.  The LSL method is to instead use the ROM-based data generated internal fields, which we denote ${\bf u}$.  As one would expect, these data generated fields would be much more accurate were we to have the full MIMO data $\{ F^{ij} \}$ for $i,j=1,\dots K$ instead of the SISO only data $\{ F^{ii}\}$. This is precisely the motivation for this work. We propose the following algorithm:
\begin{alg}[Data completion for ROM-based internal solutions] 
\label{alg:romcompute}
\\ Given monostatic data $ \{ F^{ii}(k \tau) \}$ for $k=0,\dots, 2n-2$ for sources  $i=1,\ldots, K$:
\begin{enumerate}
\item SISO step: Use LSL for true SISO data $F$ to generate approximate ${\bf q}$ and internal fields ${\bf u}^{(i)}$.
\item Lifting step: given approximate ${\bf q}$ and internal solutions ${\bf u}^{(i)}$, compute ${\bf F}$ for all inaccessible pairs $(i;j)$  \label{liftingstep}
\item MIMO step: given a full ${\bf F}$ and its mass matrix ${\bf M}$, compute new internal fields ${\bf u}^{(i)}$ and new approximate ${\bf q}$ from Lippmann-Schwinger.
\item Exit or return to lifting step \ref{liftingstep}.
\end{enumerate}
\end{alg}

\section{Step 1: The LSL algorithm for SISO data}
This step is precisely to perform the algorithm from \cite{DrMoZa3}. We include a short description here, see \cite{DrMoZa3} and also \cite{DrMoZa},\cite{DrMoZa2} for more details. 
We assume we are given only the diagonal response 
\begin{equation}\label{datadefradar} F^{jj} (k\tau) = \int_\Omega g_j(x) \cos{(\sqrt{A}k\tau)} g_j(x) dx ,\end{equation}
for $j=1,\dots, K$, $k=0,\ldots, 2n-1$.
For each source, we will compute a mass matrix and corresponding approximate internal snapshots. 
Let $u^{(j)}(x,t)$ be the true internal fields corresponding to source $g_j$, and consider the snapshots $$u^{(j)}_k = u^{(j)}(k\tau,x)$$  for $k=0,\ldots,2n-1$. The $n\times n$  mass matrix corresponding to this source is defined by
\begin{equation} \label{massmatrixj} M^{(j)}_{kl}= \int_\Omega u^{(j)}_k u^{(j)}_l  dx \end{equation}
for $k,l =0,\ldots,n-1$. Given the expression (\ref{exactsolution}) for the exact solution, we have that,
\begin{equation} \label{massmatrixj2} M^{(j)}_{kl} = \int_\Omega  g_j(x) \cos{(\sqrt{A}k\tau)} \cos{(\sqrt{A}l\tau)} g_j(x) dx.  \end{equation}
The steps are as follows. 
\begin{enumerate}
\item Obtain the mass matrices from the data: Compute
\begin{equation} \label{massfromdataj} M^{(j)}_{kl} = \frac{1}{2} \left( F^{jj}((k-l)\tau) + F^{jj}((k+l)\tau) \right).\end{equation} 
This formula follows from (\ref{datadefradar}), (\ref{massmatrixj2}),  and the cosine angle sum formula. 
\item  Orthogonalize the snapshots: Compute the Cholesky decompositions  (note each $M^{(j)}$ is positive definite),
$$M^{(j)}=(U^{(j)})^\top{U^{(j)}}$$
where each $U^{(j)}$ is upper triangular, for $j=1,\ldots, m$. 
This allows us to define the sequentially orthogonalized snapshots 
\begin{equation}\label{vdefj}  \vec{v}^{(j)} = \vec{u}^{(j)}( U^{(j)})^{-1} , \end{equation}
where $ \vec{u}^{(j)}$ is the column vector of the true (still unknown) snapshots.
\item Repeat for the background medium: 
Compute the background snapshots $\{ u^{(j),0}\} $, the mass matrices
\begin{equation} \label{massmatrix0j} M^{(j),0}_{kl}= \int_\Omega  u^{(j),0}_k  u^{(j),0}_l dx, \end{equation}
their Cholesky decompositions $$M^{(j),0}=(U^{(j),0})^\top U^{(j),0},$$
and the orthogonalized background snapshots, 
\begin{equation}\label{v0defj}  \vec{v}^{(j),0} = \vec{u}^{(j),0} (U^{(j),0})^{-1} . \end{equation}
\item Compute data generated internal snapshots: In this step we replace the unknown orthogonalized snapshots $\vec{v}^{(j)}$ with the orthogonalized background snapshots $\vec{v}^{(j),0}$. {They should be close in the sense of spherical averages; see \cite{DrMoZa3} Section 4.2 for a discussion and Figure 7.7 in Section \ref{sec:numerical} here for a numerical example comparing internal solutions.} 
From (\ref{vdefj}) and (\ref{v0defj}), we have our approximations to the internal snapshots:
  \begin{eqnarray} \label{internalj} \vec{\bf{u}}^{(j)} &=&\vec{v}^{(j),0} U^{(j)}  \\ &=& \vec{u}^{(j),0} (U^{(j),0})^{-1}U^{(j)} . \nonumber\end{eqnarray}
\item Invert the Lippmann-Schwinger equation: 
In order to obtain an estimate for the unknown $q$, we now use the Lippmann-Schwinger equation. For each time step $k\tau$, $k= 0,\dots, n-1$, and for each source $j=1,\ldots, m$, we have that 
\begin{equation} \label{Lipschj} F^{jj}_{0}(k\tau) - F^{jj}(k\tau) = \int_0^{k\tau}\int_\Omega w^{(j),0}(x,k\tau-t) u^{(j)}(x,t) q(x) dxdt .\end{equation}
The Lippmann-Schwinger-Lanczos method is to replace $u^{(j)}$ in the above by its data generated approximation $\vec{\bf{u}}^{(j)}$ from (\ref{internalj}).  Let ${\bf u}^{(j)}(x,t)$ be a time interpolant of $\vec{\bf{u}}^{(j)}$, and we solve
\begin{equation} \label{LipschLj} F^{jj}_0(k\tau) - F^{jj}(k\tau) = \int_0^{k\tau}\int_\Omega w^{(j),0}(x,k\tau-t) {\bf{u}}^{(j)}(x,t) q(x) dxdt ,\end{equation}
for $j=1,\ldots,K$ , $k= 0,\ldots ,n-1$, yielding $nm$ equations to be inverted to find ${\bf q}$, an approximation to $q$. 
\end{enumerate}

\section{Step 2: Data completion and its mass matrix}
In this step we assume we are given an approximation to the internal snapshots ${\bf u}^{(j)}(x,t)$ for $k=0,\ldots n-1$, $j=1,\ldots K$, along with an approximation $\bf{q}$ to the unknown $q$. We now use (\ref{eq:LSTD}) to approximate the off diagonals of $F$. 
\begin{enumerate}
\item Data lifting: For $i\neq j$, compute approximations ${\bf F}^{ij}$ to $F^{ij}$ 
\begin{equation} \label{liftingLip}   {\bf F}^{ij}(k\tau) := F^{ij}_{0}(k\tau) - \int_0^{k\tau}\int_\Omega w^{(j),0}(x,k\tau-t) {\bf u}^{(i)}(x,t) {\bf q}(x) dxdt. \end{equation}
Define the diagonals \begin{equation} {\bf F}^{ii} =F^{ii} \end{equation}
by the actual data values. 
What we have now is the lifted data  
$$ \{ {\bf F}^{ij}(k \tau) \} \  \mbox{ for } \  k=0,\dots, n-1 \ \  i,j=1,\ldots K.$$ 
Note that the number of time steps is reduced by half compared to the original data (or the data from the previous lifting iterate). 
\item MIMO mass matrix: Since ${\bf F}$ is now a full $K\times K$ matrix, we use the block version of (\ref{massfromdataj}) (see for example \cite{BoGaMaZi}) to compute the MIMO mass matrix 
\begin{equation} \label{massfromdatablock} {\bf M}_{kl} = \frac{1}{2} \left( {\bf F}((k-l)\tau) + {\bf F}((k+l)\tau) \right)\end{equation} 
for $k,l = 0,\ldots, \lfloor (n-1)/2 \rfloor$.
Here ${\bf M}_{kl}$ represents the $kl$th $K\times K$ block of ${\bf M}$.
\item Regularization of ${\bf M}$: 
Since ${\bf M}$ is not a true mass matrix, it will likely not be symmetric positive definite, so regularization is required. We regularize by taking the symmetric part and thresholding: all eigenvalues below a chosen small  $\epsilon_0> 0 $ are set to $\epsilon_0$. We note the result here is not very sensitive to the choice of $\epsilon_0$,  provided that $\epsilon_0$ is smaller than the minimum positive eigenvalue $\lambda^+_{min}$ of ${\bf M}$. In our experiments we pick $\epsilon_0=\sqrt{10^{-12}\lambda^+_{max}\lambda^+_{min}}$  where $\lambda^+_{max}$ is the maximum positive eigenvalue of ${\bf M}$. We note that this $\epsilon_0$ corresponds to the geometric mean of $\lambda^+_{min}$ and the minimum value $10^{-12}\lambda^+_{max}$ for stable linear algebraic operations with modified ${\bf M}$ in machine double precision. We abuse notation and continue to denote this regularized mass matrix by ${\bf M}$. 
\end{enumerate}

\section{Step 3: The MIMO LSL algorithm for completed data} 
In this step we assume we have a positive definite mass matrix ${\bf M}$ computed from a full ${\bf F}$, where ${\bf M}$ is corresponding to time steps $k=0,\ldots,N$. The first time that we perform this step we will have $N=  \lfloor (n-1)/2 \rfloor$; subsequent steps will be halved in size.
\begin{enumerate}
\item Orthogonalize  the snapshots: We do a block Cholesky  (with $K\times K$ blocks) decomposition, 
$$ {\bf M} = U^\top U$$ 
where $U$ is block upper triangular.  We note that there is some ambiguity in the choice of the blocks; we choose them so that all resulting orthogonalized functions 
$$ \vec{v} = \vec{u} U^{-1} .$$
are all mutually orthogonal. 
\item Repeat for the background medium: Compute the background MIMO mass matrix corresponding to the same time steps, its Cholesky decomposition
$$ M^0 = (U^0)^\top U^0,$$ and the orthogonalized background snapshots
$$ \vec{v}^0 = \vec{u}^0 (U^0)^{-1}. $$
\item Compute new data generated internal snapshots:   
$$\vec{\bf{u}} = \vec{u}^0 (U^0)^{-1} U .$$
An important point is that these internal fields have been improved substantially compared to those generated in Step 1. 
\item Invert the Lippman-Schwinger equation: 
Let ${\bf u}^{(j)}(x,t)$ be a time interpolant of the fields in $\vec{\bf{u}}$ corresponding to the $j$th source, and solve the Lippmann-Schwinger equation 
\begin{equation} \label{LipschLjcompleted} F^{jj}_0(k\tau) - F^{jj}(k\tau) = \int_0^{k\tau}\int_\Omega w^{(j),0}(x,k\tau-t) {\bf{u}}^{(j)}(x,t) q(x) dxdt \end{equation}
to obtain a new approximation ${\bf q}$ to the unknown $q$.  In this Lippmann-Schwinger step we use only the true SISO data, since our lifted off diagonals do not add any new information; they were used only to improve the internal fields. 
\end{enumerate}

\section{Step 4: loop for iterations}
We note that our experiments suggest that whether or not we need to proceed with iterations and return to step \ref{liftingstep} depends on the number of targets to resolve. Our observations show that the number of iterations should be one less than the maximum number of objects with internal multiples. In particular, for the problems with two targets we will consider section \ref{sec:numerical}, no loop is needed. However, there is benefit from having an extra iteration for the three-objects problem (see our last numerical example). We can vaguely explain this as follows: the Born solution (as well as LSL without completion) successfully resolve the closest target only, while the images of deeper ones are polluted by multiple echoes. Then, every iteration of LSL with completion cleans up these multiples from the deeper and deeper objects one by one. In this sense there is a similarity between our approach and the Inverse Born Series \cite{IBS}, however, a more precise connection needs to be explored further.

\section{Numerical Experiments}
\label{sec:numerical}

\begin{figure}[htb]
\centering
\includegraphics[width=.55\textwidth]{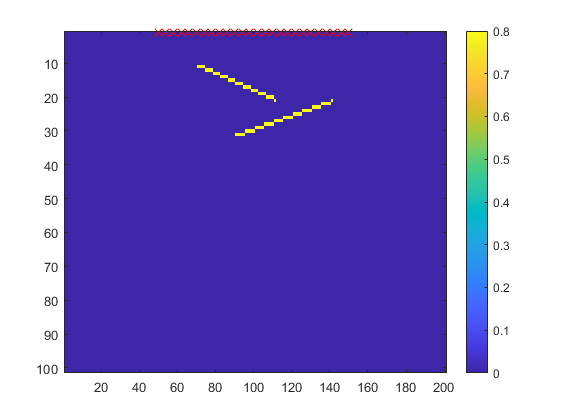}
\vspace{-2ex}
\caption{True model for imaging of two targets in homogeneous background. Red crosses show the locations of SAR sources/receivers.}
\label{fig:med}
\end{figure}

\begin{figure}[htb]
\includegraphics[width=0.9\textwidth]{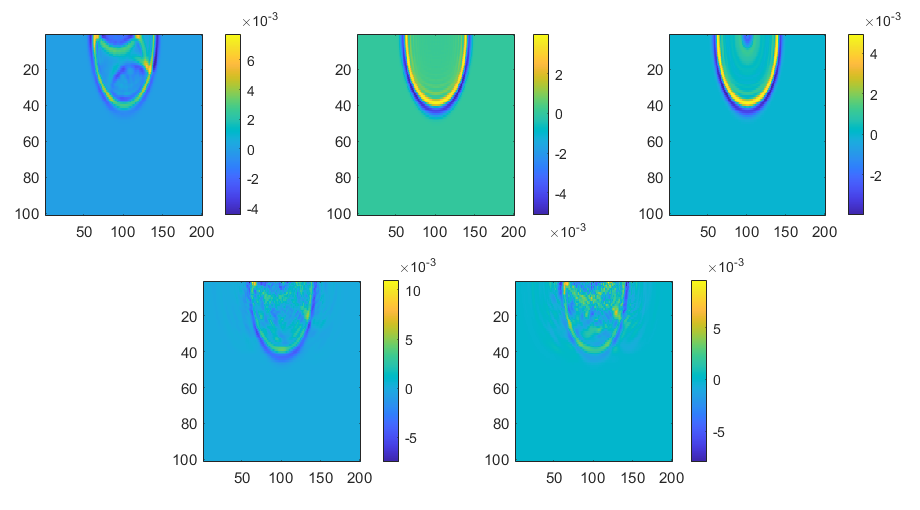}
\vspace{-2ex}
\caption{  
While the background solution (top middle) doesn't capture any reflections of the true solution (top left), the LSL solution without data completion captures them in spherical averages only (top right). The data completion step improved the internal solution reconstruction (bottom left) significantly. The second iteration (bottom right) sharpened the profile, though all of the reflections were already captured during the first iteration.}
\label{fig:intsol_compl}
\end{figure}

\begin{figure}[htb]
	\centering
	\begin{tabular}{cc}
	True & Born \\
		\includegraphics[width=0.39\textwidth]{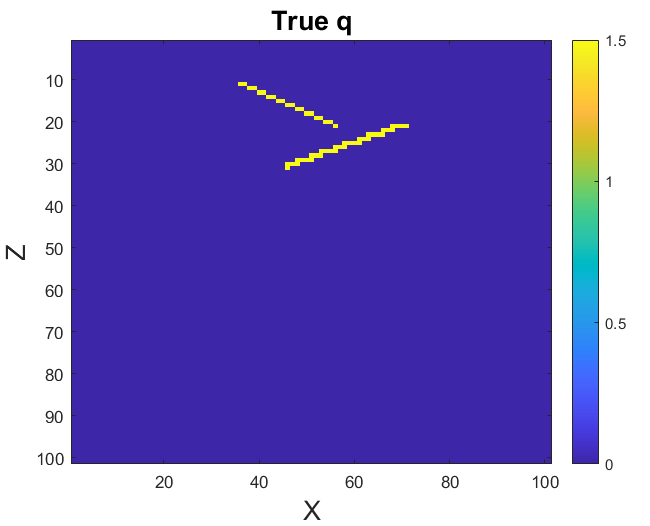} &
		\includegraphics[width=0.4\textwidth]{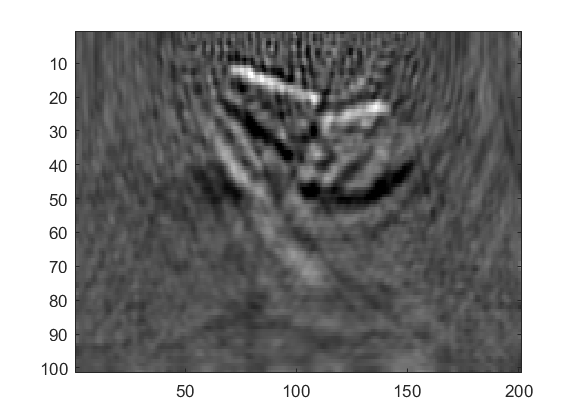}\\
		LSL w/o completion & LSL with completion \\
		\includegraphics[width=0.4\textwidth]{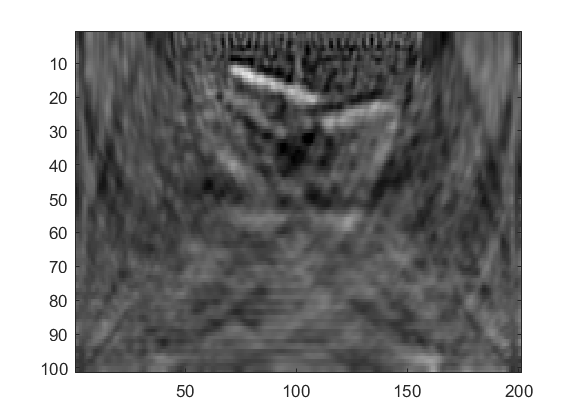} &
		\includegraphics[width=0.4\textwidth]{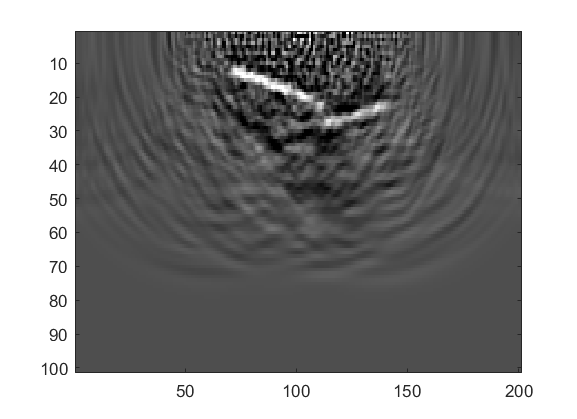}\\
	\end{tabular}
	\begin{tabular}{c}
	2nd iteration of data completion fo LSL\\
		\includegraphics[width=0.39\textwidth]{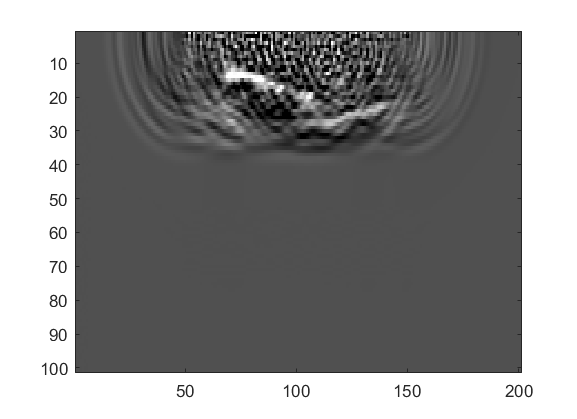}
	\end{tabular}
	\caption{2D example. LSL solution without data completion (bottom left) improves upon the Born solution (top right) . One step of LSL data completion (bottom right) allowed us to improve the image significantly. The second iteration didn't improve the image for this simple model.}
	\label{fig:numex1}
\end{figure}

\begin{figure}[htb]
	\centering
	\begin{tabular}{cc}
	True & Born \\
		\includegraphics[width=0.4\textwidth]{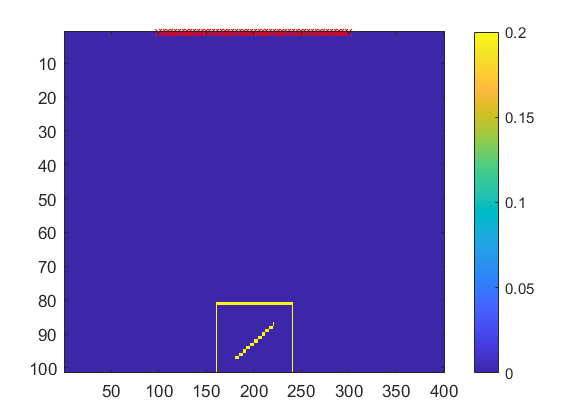} &
		\includegraphics[width=0.4\textwidth]{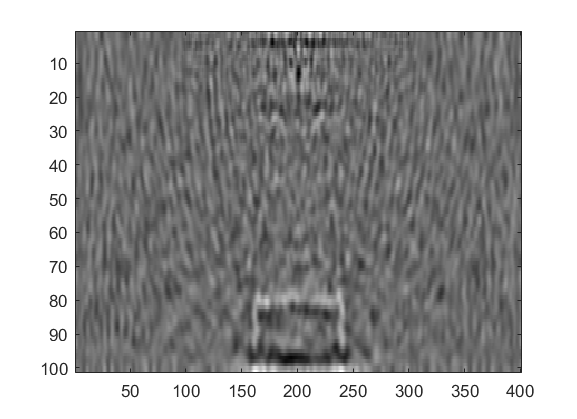}\\
		LSL w/o completion & LSL with completion \\
		\includegraphics[width=0.4\textwidth]{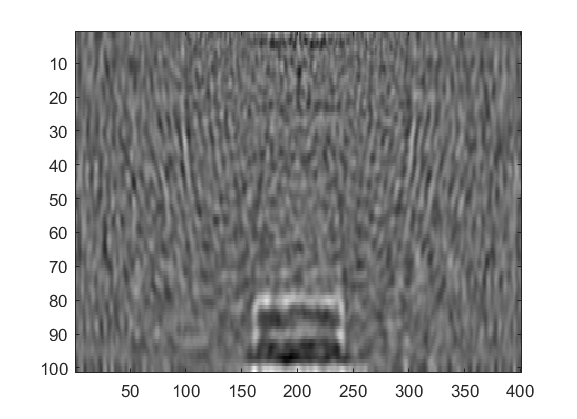} &
		\includegraphics[width=0.4\textwidth]{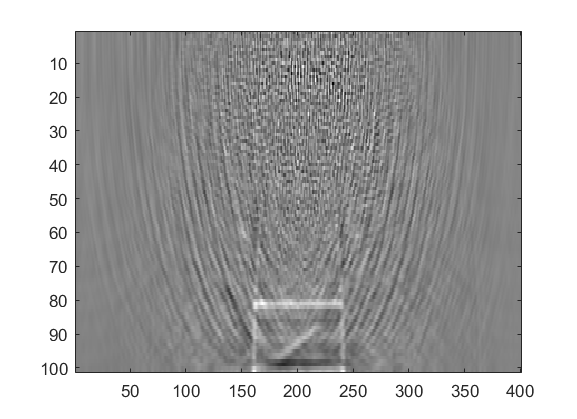}\\
	\end{tabular}
	\begin{tabular}{c}
		LSL with completion for noisy data\\
		\includegraphics[width=0.4\textwidth]{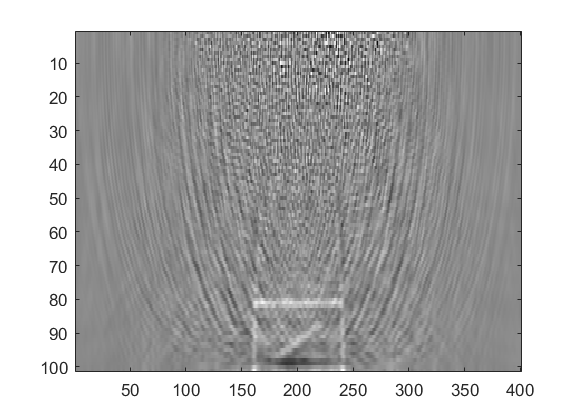}\\
	\end{tabular}
	\caption{2D example: Red crosses in the true model (top left) show the locations of SAR sources/receivers. LSL solution without data completion (middle left) improves upon the Born solution (top right) . One step of LSL data completion (middle right) allowed us to resolve the hidden object. Adding 5\% noise didn't worsen the image (bottom) significantly.}
	\label{fig:numexbox}
	
\end{figure}

\begin{figure}[htb]
	\centering
	\begin{tabular}{cc}
	True & Born \\
		\includegraphics[width=0.4\textwidth]{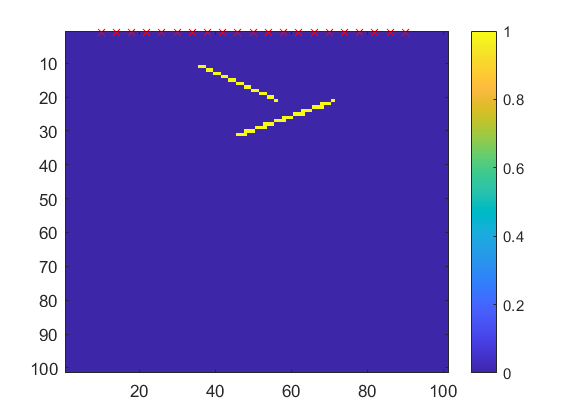} &
		\includegraphics[width=0.4\textwidth]{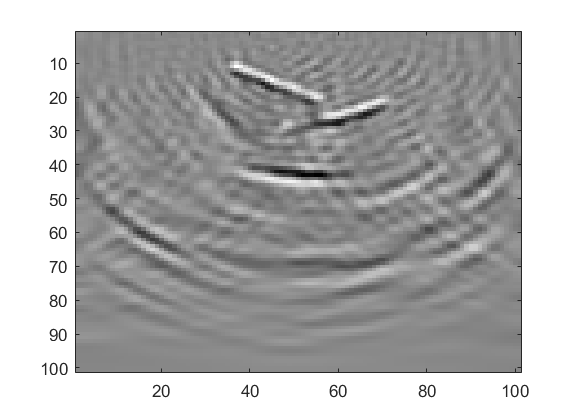}\\
		LSL w/o completion & LSL with completion \\
		\includegraphics[width=0.4\textwidth]{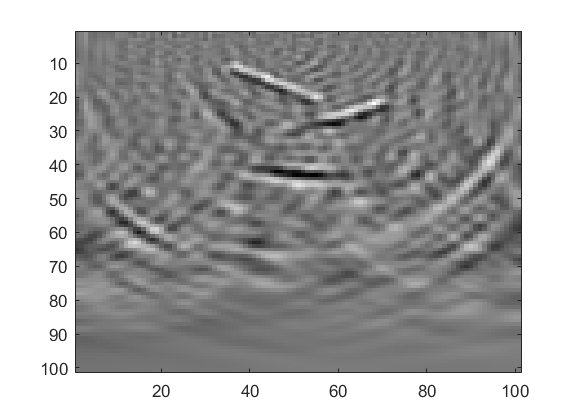} &
		\includegraphics[width=0.4\textwidth]{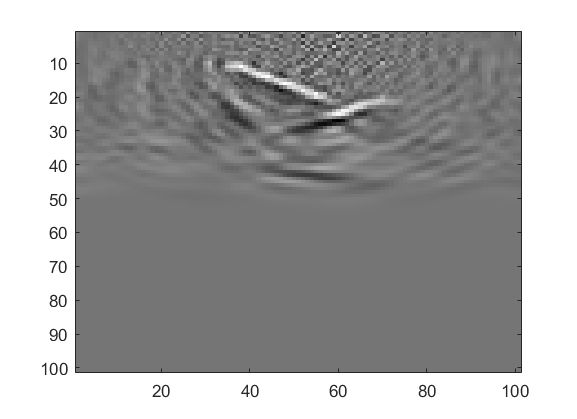}\\
	\end{tabular}
	\caption{ 2.5D example: Red crosses in the true model (top left) show the locations of SAR sources/receivers. LSL solution without data completion (bottom left) , Born solution (top right). One step of LSL data completion (bottom right) yields a dramatic improvement in the image. }
	\label{fig:numex2}
\end{figure}

In our first numerical experiment, we considered the imaging of two elongated targets in a homogeneous background (see Fig. \ref{fig:med}). To obtain synthetic data, we discretized \ref{waveeq} in $\Omega=[0; 200]\times [0;100]$ using finite-differences on a $400\times 200$ 2D grid. The obtained discrete problem was solved for 27 positions of radar that were emitting a non-modulated Gaussian pulse (radar positions are marked by crosses in Fig. \ref{fig:med}). To solve the Lippmann-Schwinger equation, we discretized it using quadrature on a $200\times 100$ grid. First, we examine how data completion improves the accuracy of the internal fields. In Fig. \ref{fig:intsol_compl} we see the true internal solution (top left) and the background solution (top middle), which is used for the Born approximation in (\ref{eq:LSTD}). The ROM-based SISO internal solution without data completion is shown in the top right.  While the background solution doesn't have any reflections, the internal solution we obtained without data completion captures the reflections in the sense of spherical averages. The internal solution computed after data completion is shown in the bottom left. As one can observe, MIMO data completion improves the accuracy of the internal solution reconstruction significantly. We also note that the second iteration made the profile sharper (see bottom right), however, all of the multiples were already captured during the first iteration. 

In Figure \ref{fig:numex1} we compare reconstructions for this first example, where we see the Born image (top right) and LSL images without and with data completion (middle left and right respectively). We note that the ROM step had to be performed in a structure-preserving way in order to obtain a stable internal solution. The data completion for this first example results in a significantly sharper image of two targets, with almost no ghost reflections. As was expected, the second iteration (bottom plot) didn't yield any improvement for this two target model.   
{We note that in each of the reconstructions we chose the SVD truncation level which was best for each separately, which is why the frequency of the artifacts appears higher for the LSL with completion. In all cases the artifacts seen in these pictures could possibly be further diminished by using sparse or other regularization techniques.}

For our second example we considered another two dimensional target (see top left plot in Figure \ref{fig:numexbox}) with 60 sources and receivers located at the top and emitting a non-modulated Gaussian pulse. We used an $800\times 300$ finite-difference grid for discretization of the forward problem in $\Omega=[0;400]\times [0;150]$ to generate synthetic data. Then, data-driven internal solutions were obtained and we solved the Lippmann--Schwinger equation using the quadrature on a $400\times 150$ grid. Similar to the previous example, the data completion (middle right) really sharpens the image compared to Born (top right) and LSL without completion (middle left), and in this case allows us to resolve the image inside of the box. In this experiment we added the result with 5\% relative random noise (see bottom plot) to show a low sensitivity of the approach to the measurement errors. 

For our next example, we considered a 2.5D problem, where our waves are fully three dimensional but the medium is assumed to be homogeneous in the transversal direction (see top left in Fig. \ref{fig:numex2}). The forward problem was discretized on a $200\times 200\times 200$ grid in a 3D domain $\Omega=[0;100]^3$. Similar to the first example, we used 27 sources/receivers with a non-modulated Gaussian pulse and solved the discrete problems for them to obtain synthetic data. We then computed the ROM-based approximants of the internal solutions and plugged them in into the Lippmann-Schwinger equation, which was then solved for a 2D image via quadratures on a $100\times 100$ grid. We see the reconstructions in Figure \ref{fig:numex2}. Born (top right) and LSL without completion (bottom left) both show significant ghost reflections, but these are greatly reduced with one step of LSL with data completion (bottom right).

\begin{figure}[htb]
\centering
\includegraphics[width=.55\textwidth]{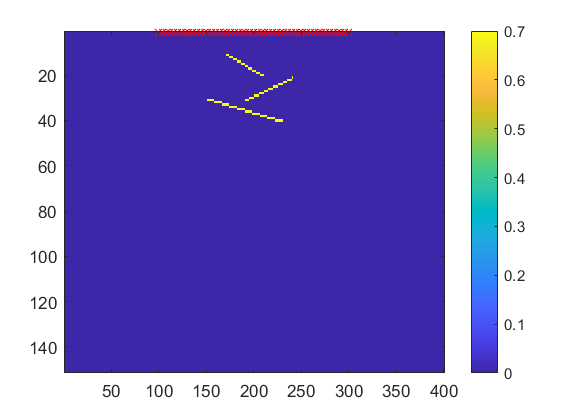}
\vspace{-2ex}
\caption{True model for imaging of three targets in homogeneous background. Red crosses show the locations of SAR sources/receivers. }
\label{fig:med_3obj}
\end{figure}

\begin{figure}[htb]
\centering
\includegraphics[width=0.9\textwidth]{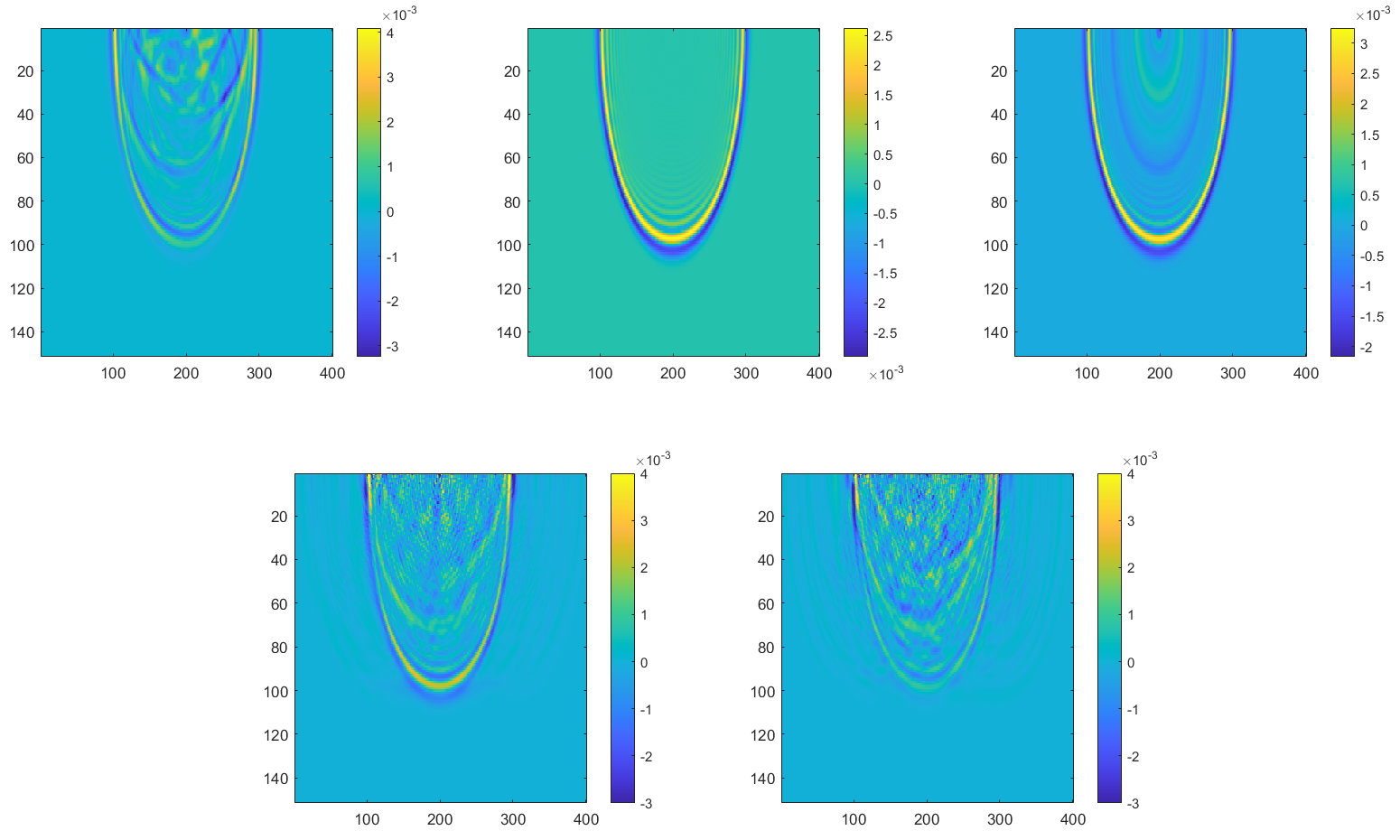}
\vspace{-2ex}
\caption{  
Similar to the problem with 2 objects, the background solution (top middle) doesn't capture any reflections of the true solution (top left) while the LSL solution without data completion captures them in spherical averages only (top right). The data completion step improved the internal solution reconstruction (bottom left) significantly, and the second iteration corrected for the missing internal multiples even further (bottom right).}
\label{fig:intsol_compl_3obj}
\end{figure}

\begin{figure}[htb]
	\centering
	\begin{tabular}{cc}
	Born &LSL w/o completion\\
		\includegraphics[width=0.4\textwidth]{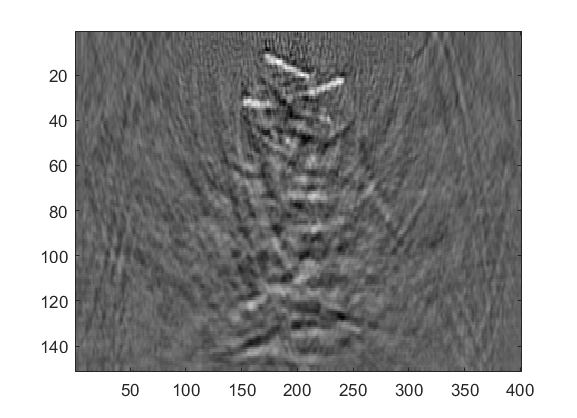} &
		\includegraphics[width=0.4\textwidth]{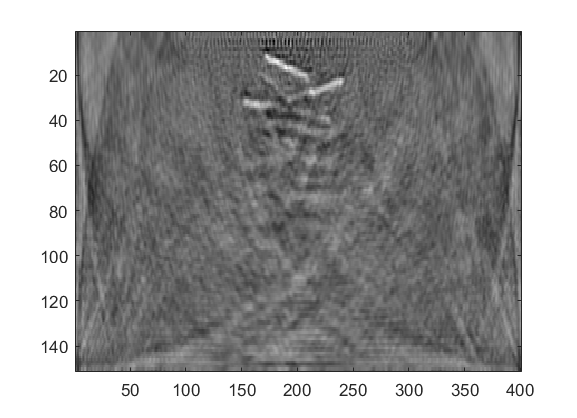}\\
		1 iteration of LSL with completion &  2 iterations of LSL with completion\\
		\includegraphics[width=0.4\textwidth]{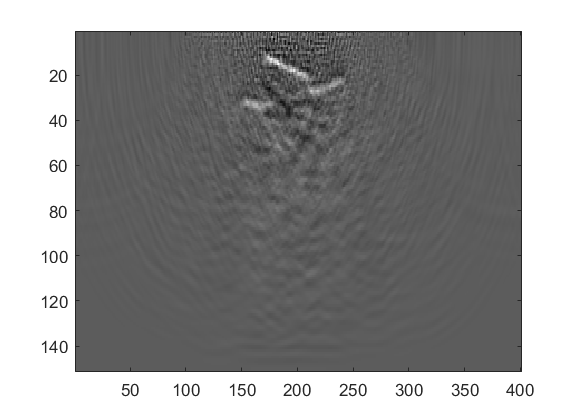} &
		\includegraphics[width=0.4\textwidth]{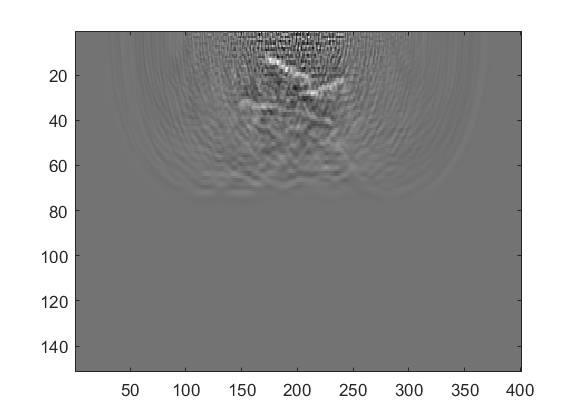}\\
	\end{tabular}
	\caption{2D example with 3 objects: LSL solution without data completion (top right) improves upon the Born solution (top left), however the image is still polluted by multiple echoes. One step of LSL data completion (bottom left) allowed us to improve the image significantly, but the deepest object only became visible after 1 extra iteration of data completion in LSL (bottom right). }
	\label{fig:3obj}
\end{figure}

For all of the examples above, adding an extra iteration of data completion in LSL will not improve the image. Indeed, as we noted before, iterations for LSL with data completion are beneficial when there are more than two objects that produce internal multiples. To illustrate this, we considered the 2D problem with 3 objects shown in Fig. \ref{fig:med_3obj}. Our modeling and inversion setups were precisely the same as in the second example, i.e. an $800\times 300$ grid in $\Omega=[0;400]\times [0;150]$, non-modulated Gaussian pulses emitted by 60 sources/receivers and quadrature for the Lippmann-Schwinger equation on a $400\times 150$ grid. In Fig. \ref{fig:intsol_compl_3obj} we show how the extra iteration (bottom right) of LSL with data completion improved the internal solution (bottom left), though both are significantly closer to the exact solution (top left) than Born (top middle) and LSL without completion (top right). This improved internal solution helped to reveal the deepest object in the model (see bottom right plot in Fig.\ref{fig:3obj}), which was neither visible in the Born image (top left), nor in the LSL images without completion (top right) and with one iteration of completion (bottom left). 

\section{Conclusions}

In the ROM based Lippmann-Schwinger-Lanczos method, one constructs data driven internal solutions, and uses them in the Lippmann-Schwinger integral. When the given data is monostatic, as is the case for SAR, these internal solutions are only accurate up to spherical averages. We saw here that adding the step of ROM based data completion can dramatically improve the accuracy of the internal solutions, leading to sharpened reconstructions. Moreover, additional iterations of data completion in LSL can improve images further in the case of complex models with more than two objects producing internal multiples.

Many practically important problems have more complicated data matrix patterns than shown in (\ref{fullMIMO}-\ref{SAR}). For example, in the case of tomographic ultrasound, the data acquisition may be block diagonal. The method we described in this work can be extended many of those cases as well \cite{Roberts:1980:LMR}.  
{The LSL approach can also be extended in some cases to variable wave speed problems, in particular for small slow speed perturbations \cite{BaChDrMoZa}, or more generally by combining with iteration \cite{BoGaMaZi}. }
We plan to investigate extension of  this data-completion approach to these and more general setups in future work.

\thanks{{\bf Acknowledgments} 
The authors are  grateful to Liliana Borcea, Alex Mamonov and J\"orn Zimmerling for productive discussions that inspired this research.
V. Druskin was partially supported by AFOSR grants FA 955020-1-0079, FA9550-20-1-0079,  and NSF grant  DMS-2110773. M. Zaslavsky was partially supported by  AFOSR grant  FA9550-20-1-0079. S. Moskow was partially supported by NSF grants DMS-2008441 and DMS-2308200.

\bibliography{biblio,biblio6,functions-of-matrices,graphbib6,galerkincitations}

\def\noopsort#1{}\def\l{\char32l}\def\v#1{{\accent20 #1}}
  \let\^^_=\v\def\hbk{hardback}\def\pbk{paperback}
\begin{thebibliography}{10}

\bibitem{BaChDrMoZa}
J.~Baker, E.~Cherkaev, V.~Druskin, S.~Moskow, and M.~Zaslavsky.
\newblock Doubly regularized lippmann-schwinger-lanczos algorithm for inverse
  scattering problems in the frequency domain.
\newblock submitted.

\bibitem{borcea2020reduced}
L.~Borcea, V.~Druskin, A.~Mamonov, M.~Zaslavsky, and J.~Zimmerling.
\newblock Reduced order model approach to inverse scattering.
\newblock {\em SIAM Journal on Imaging Sciences}, 13(2):685--723, 2020.

\bibitem{borcea2011resistor}
Liliana Borcea, Vladimir Druskin, Fernando Guevara~Vasquez, and Alexander~V.
  Mamonov.
\newblock Resistor network approaches to electrical impedance tomography.
\newblock {\em Inverse Problems and Applications: Inside Out II, Math. Sci.
  Res. Inst. Publ}, 60:55--118, 2011.

\bibitem{BoDrMaMoZa}
Liliana Borcea, Vladimir Druskin, Alexander~V. Mamonov, Shari Moskow, and
  Mikhail Zaslavsky.
\newblock Reduced order models for spectral domain inversion: embedding into
  the continuous problem and generation of internal data.
\newblock {\em Inverse Problems}, 36(5), 2020.

\bibitem{borcea2014model}
Liliana Borcea, Vladimir Druskin, Alexander~V. Mamonov, and Mikhail Zaslavsky.
\newblock A model reduction approach to numerical inversion for a parabolic
  partial differential equation.
\newblock {\em Inverse Problems}, 30(12):125011, 2014.

\bibitem{borcea2017untangling}
Liliana Borcea, Vladimir Druskin, Alexander~V Mamonov, and Mikhail Zaslavsky.
\newblock Untangling nonlinearity in inverse scattering with data-driven
  reduced order models.
\newblock {\em Inverse Problems}, 34(6):065008, 2018.

\bibitem{borcea2019robust}
Liliana Borcea, Vladimir Druskin, Alexander~V. Mamonov, and Mikhail Zaslavsky.
\newblock Robust nonlinear processing of active array data in inverse
  scattering via truncated reduced order models.
\newblock {\em Journal of Computational Physics}, 381:1--26, 2019.

\bibitem{BoGaMaZi}
Liliana Borcea, Josselin Garnier, Alexander~V. Mamonov, and J\"{o}rn
  Zimmerling.
\newblock Waveform inversion with a data driven estimate of the internal wave.
\newblock {\em SIAM Journal on Imaging Sciences}, 16(1):280--312, 2023.

\bibitem{Borcea2021ReducedOM}
Liliana Borcea, Josselin Garnier, Alexander~V. Mamonov, and J{\"o}rn~T.
  Zimmerling.
\newblock Reduced order model approach for imaging with waves.
\newblock {\em ArXiv}, abs/2108.01609, 2021.

\bibitem{Burfeindt2}
Matthew~J. Burfeindt, Hatim~F. Alqadah, and S.~Ziegler.
\newblock Experimental phase-encoded linear sampling method imaging with a
  single transmitter and receiver.
\newblock In {\em IEEE Open Journal of Antennas and Propagation}, pages 1--17,
  2024.

\bibitem{CaCoMo}
Fioralba Cakoni, David Colton, and Peter Monk.
\newblock {\em The Linear Sampling Method in Inverse Electromagnetic
  Scattering}.
\newblock Society for Industrial and Applied Mathematics, 2011.

\bibitem{cheney2021procedure}
Margaret Cheney, Matthew~J Burfeindt, and NAVAL RESEARCH LAB WASHINGTON
  DCColorado~State University.
\newblock A procedure for suppressing multiple scattering.
\newblock 2021.

\bibitem{druskin2018nonlinear}
V.~Druskin, A.V. Mamonov, and M.~Zaslavsky.
\newblock A nonlinear method for imaging with acoustic waves via reduced order
  model backprojection.
\newblock {\em SIAM Journal on Imaging Sciences}, 11(1):164--196, 2018.

\bibitem{druskin2016direct}
Vladimir Druskin, Alexander~V. Mamonov, Andrew~E. Thaler, and Mikhail
  Zaslavsky.
\newblock Direct, nonlinear inversion algorithm for hyperbolic problems via
  projection-based model reduction.
\newblock {\em SIAM Journal on Imaging Sciences}, 9(2):684--747, 2016.

\bibitem{DrMoZa}
Vladimir Druskin, Shari Moskow, and Mikhail Zaslavsky.
\newblock Lippmann-schwinger-lanczos algorithm for inverse scattering problems.
\newblock {\em Inverse Problems}, apr 2021.

\bibitem{DrMoZa2}
Vladimir Druskin, Shari Moskow, and Mikhail Zaslavsky.
\newblock On extension of the data driven rom inverse scattering framework to
  partially nonreciprocal arrays.
\newblock {\em Inverse Problems}, 8(38), 2022.

\bibitem{DrMoZa3}
Vladimir Druskin, Shari Moskow, and Mikhail Zaslavsky.
\newblock Reduced order modeling inversion of monostatic data in a
  multi-scattering environment.
\newblock {\em SIAM Journal on Imaging Sciences}, 17(1), 2024.

\bibitem{gilman2015mathematical}
Mikhail Gilman and Semyon Tsynkov.
\newblock A mathematical model for sar imaging beyond the first born
  approximation.
\newblock {\em SIAM Journal on Imaging Sciences}, 8(1):186--225, 2015.

\bibitem{Kepley2016GeneratingVI}
Paul Kepley, Lauri Oksanen, and Maarten~V. de~Hoop.
\newblock Generating virtual interior point source traveltimes and redatuming
  using boundary control.
\newblock {\em Seg Technical Program Expanded Abstracts}, 2016.

\bibitem{10352835}
Jerry~T. Kim, Eric~L. Mokole, and Margaret Cheney.
\newblock Concept of iterative time-reversal radar (itrr).
\newblock In {\em 2023 IEEE Conference on Antenna Measurements and Applications
  (CAMA)}, pages 411--416, 2023.

\bibitem{9357476}
Michael~V. Klibanov, Alexey~V. Smirnov, Vo~A. Khoa, Anders~J. Sullivan, and
  Lam~H. Nguyen.
\newblock Through-the-wall nonlinear sar imaging.
\newblock {\em IEEE Transactions on Geoscience and Remote Sensing},
  59(9):7475--7486, 2021.

\bibitem{lagergren2021deep}
John Lagergren, Kevin Flores, Mikhail Gilman, and Semyon Tsynkov.
\newblock Deep learning approach to the detection of scattering delay in radar
  images.
\newblock {\em Journal of Statistical Theory and Practice}, 15(1):14, 2021.

\bibitem{Malcolm2007IdentificationOI}
Alison~E. Malcolm, Maarten~V. de~Hoop, and Henri Calandra.
\newblock Identification of image artifacts from internal multiples.
\newblock {\em Geophysics}, 72, 2007.

\bibitem{IBS}
Shari Moskow and John~C. Schotland.
\newblock {\em 12. Inverse Born series. The Radon Transform: The First 100
  Years and Beyond, edited by Ronny Ramlau and Otmar Scherzer}, pages 273--296.
\newblock De Gruyter, Berlin, Boston, 2019.

\bibitem{Roberts:1980:LMR}
J.~D. Roberts.
\newblock Linear model reduction and solution of the algebraic {Riccati}
  equation by use of the sign function.
\newblock {\em Internat. J. Control}, 32(4):677--687, 1980.
\newblock First issued as report CUED/B-Control/TR13, Department of
  Engineering, University of Cambridge, 1971.

\bibitem{Virieux2016AnIT}
Jean Virieux.
\newblock An introduction to full waveform inversion.
\newblock 2016.

\bibitem{WaAnLi}
Francis Watson, Daniel Andre, and William R~B Lionheart.
\newblock Resolving full-wave through-wall transmission effects in multi-static
  synthetic aperture radar.
\newblock {\em Inverse Problems}, 2024.

\end{thebibliography}
\end{document}